\documentclass[11pt]{paper}

\usepackage{amssymb,amsmath,enumerate,theorem,epsfig,rotating,graphics,changebar,eepic}
\usepackage{amsfonts}
\usepackage{a4,latexsym,parskip}
\usepackage{hyperref}
\hypersetup{pdfauthor=MS} \hypersetup{pdftitle=Quintic}

\newcommand{\PP}{\mathbb{P}}
\newcommand{\C}{\mathbb{C}}
\newcommand{\Q} {\mathbb{Q}}

\newcommand{\F}{\mathbb{F}}
\newcommand{\Z}{\mathbb{Z}}
\newcommand{\p}{\mathfrak{p}}

\newcommand{\qed}{\hfill \ensuremath{\Box}}
\newcommand{\NS}{\mathop{\rm NS}\nolimits}
\newcommand{\A}{\mathfrak{A}}
\newcommand{\TT}{\mathfrak{T}}

\newcommand{\Frob}{\mathop{\rm Frob}\nolimits}

\theoremstyle{break} \newtheorem{Theorem}{Theorem}
\newtheorem{Proposition}[Theorem]{Proposition}
\newtheorem{Lemma}[Theorem]{Lemma}
\newtheorem{Example}[Theorem]{Example}

\newtheorem{Remark}[Theorem]{Remark}

\begin{document}
\setlength{\unitlength}{1cm}

\title{Quintic surfaces with maximum and other Picard numbers} 

\author{Matthias Sch\"utt}


\date{\today}

\maketitle

\abstract{This paper investigates the Picard numbers of quintic surfaces.
We give the first example of a complex quintic surface in $\PP^3$ with maximum Picard  number $\rho=45$. 
We also investigate its arithmetic and determine the zeta function.
Similar techniques are applied to produce quintic surfaces with several other Picard numbers that have not been achieved before.}



\section{Introduction}

This paper concerns the problem of exhibiting complex algebraic surfaces of general type with given Picard number. 
In general, there are only a few Picard numbers known to be attained within a fixed class of algebraic surfaces.
In particular it is unclear whether every Picard number satisfying Lefschetz' bound
\begin{eqnarray}\label{eq:Lef}
\rho(X)\leq h^{1,1}(X)
\end{eqnarray}
might be attained.
In this paper we concentrate on the case of quintic surfaces in $\PP^3$.
The non-trivial Hodge numbers of a quintic surface $X$ are
\[
h^{2,0}(X)=4,\;\;\; h^{1,1}(X)=45,\;\;\; h^{0,2}(X)=4.
\]
We will extend the known results greatly by providing specific examples in Section \ref{s:rho}.
Special emphasis is put on the case of maximum Picard  number.
A smooth compact complex surface $X$ is said to have \emph{maximum Picard number}
if its Picard number $\rho(X)$ attains the Lefschetz bound \eqref{eq:Lef}.
This property is a birational invariant of $X$,
and we often employ the same terminology for irreducible singular surfaces 
by considering their desingularisations.

There are a few deformation classes of smooth surfaces 
in which we know the existence of surfaces with maximum Picard number. 
Typical examples are:
surfaces with $h^{2,0}(X)=0$ (trivial case);
abelian surfaces and K3 surfaces (by the Torelli theorem);
 and certain double covers 
of rational surfaces (Persson \cite{Persson}).
In general, however, it is a question widely open
whether a given deformation class of surfaces contains a member with maximum Picard  number or not.
For instance, it has not been known whether a surface of degree $d$ in $\PP^3$
can have maximum Picard number, except for the cases $d\leq 4$ or $d=6$ (Beauville).
In this note, we address the problem when $X$ is a quintic surface in $\PP^3$,
answering a question raised by Shioda in \cite{Sh-PicV}.

\begin{Theorem}\label{thm}
The surface $Y\subset\PP^3$ defined by the equation
\[
yzw^3+xyz^3+wxy^3+zwx^3=0
\]
has exactly four $A_9$ singularities at the points where three coordinates vanish simultaneously.
Its minimal resolution $X$ has maximum Picard  number $\rho(X)=45$.
\end{Theorem}

We give below three proofs, each being of its own independent interest. 
The first proof exploits the fact that $X$ is the Galois quotient of a Fermat surface, 
thus closely following an idea of Shioda. 
For the second proof, we exhibit rational curves on $X$ which generate the N\'eron-Severi group $\NS(X)$ up to finite index.
The third proof uses the cyclic group of order $15$ acting on $X$ to show 
that the $\Q$-transcendental cycles form a one-dimensional vector space over the cyclotomic field $\Q(\zeta_{15})$. 

If we combine the first and second proofs,
we can compute the zeta function of $X$.
Meanwhile the method of the first proof allows us 
to produce quintic surfaces with intermediate Picard numbers (see Section \ref{s:rho}):

\begin{Theorem}
\label{thm2}
If $r=1,5,13$ or an odd integer between $17$ and $45$,
then there exists a quintic surface $X$ with $\rho(X)=r$.
\end{Theorem}

\section{Picard numbers of algebraic surfaces}
\label{s:over}

In this section, we review what seems to be known about Picard numbers of algebraic surfaces, especially about maximum Picard  number. In general, it is very difficult to determine the Picard number of a given surface $X$.
This problem admits several approaches that can sometimes also be combined.

Obviously, exhibiting algebraically independent divisor classes in $\NS(X)$ will give a lower bound for $\rho(X)$.
This is often achieved by computing intersection numbers and the rank of the resulting Gram matrix.
There is a trivial case where this lower bound determines $\rho(X)$:
in the case of maximum Picard  number where 
the lower bound coincides with the upper bound given by \eqref{eq:Lef} over $\C$ and by $b_2(X)$ in positive characteristic (due to Igusa).
This might serve as a first indication why the property of maximum Picard  number is so special.
In the presence of automorphisms acting non-trivially on the two-forms, these bounds have been improved by Shioda in \cite{Sh-PicV}.
For instance,
he proved that a surface $X\subset \PP^3$ of prime degree $d$,
given by an equation 
\begin{eqnarray}
\label{eq:Shioda}
w^d = f(x,y,z),
\end{eqnarray}
has Picard number $\rho(X)\leq h^{1,1}(X)-p_g(X)$.

An upper bound for the Picard number can also be obtained from specialisation.
For instance, we can start with a surface $X$ over some number field and then consider  its smooth reduction modulo some prime $\p$. 
Then $\rho(X\otimes\bar\Q)\leq \rho(X\otimes\bar\F_\p)$, 
and the latter number is bounded by the number of certain roots of the characteristic polynomial of $\Frob_\p^*$ on the \'etale cohomology groups $H^2(\bar X)$. 
At least in principle, the characteristic polynomial can be computed via Lefschetz' fixed point formula by counting points over sufficiently many extensions of $\F_\p$,
thus yielding an upper bound for both $\rho(X\otimes\bar\Q)$ and $\rho(X\otimes\bar\F_\p)$. 
The Tate conjecture predicts that this upper bound gives in fact an equality with the latter Picard number \cite{Tate}.

There is one subtlety when comparing upper and lower bound:
the parity of $b_2(X)$ prescribes  the parity of the upper bound.
For instance, smooth quintics over finite fields ought to have odd geometric Picard number by the Tate conjecture.
Along the same lines, one has even geometric Picard number for K3 surfaces over finite fields.
This complicates the search for surfaces with the opposite parity substantially.
As an illustration, consider the K3 case.
Terasoma proved as part of a more general result for complete intersections 
that there is a quartic surface in $\PP^3$ defined over $\Q$ that has Picard number one \cite{T}.
However, it took another twenty years to actually exhibit such a K3 surface explicitly in \cite{vL}.

There is one other non-trivial case where the Picard number of a surface can be computed in an intrinsic manner:
for Delsarte surfaces, one can argue with the covering Fermat surfaces by a method pioneered by Shioda \cite{Sh-Pic}.
This technique will feature prominently in this paper.
We will explain it in Section \ref{s:Del}.
In Section \ref{s:rho}, it will be used extensively to exhibit quintic surfaces with a plentitude of Picard numbers.

We shall now discuss the problem of maximum Picard  number in more detail.
The main reference is Persson's paper \cite{Persson} which established the existence for certain double covers. We will also comment on related arithmetic issues.

There is one kind of surfaces where the question of the Picard number has a trivial answer since every surface has maximum Picard  number. Recall that Lefschetz' bound (\ref{eq:Lef}) is a consequence of the more precise result that
\[
\mbox{Pic}(X) = H^2(X,\Z) \cap H^{1,1}(X).
\]
Hence $h^{2,0}(X)=0$ implies $\rho(X)=h^{1,1}(X)$. Thus we are led to consider surfaces with $h^{2,0}(X)\neq 0$. 

The problem of maximum Picard  number was classically solved 
for complex abelian surfaces and K3 surfaces by the Torelli theorem: 
Here the surfaces with maximum Picard  number are often called singular and lie dense in the moduli space. The terminology does not refer to non-smoothness, but to the surfaces being exceptional.
It is borrowed from the theory of elliptic curves with complex multiplication (CM), i.e.~with extra endomorphisms. 
In fact, there is a direct connection that gives rise to many arithmetic applications.
For details, see \cite{S-NS}, \cite{SI}, \cite{SM}.
In this spirit, we will also investigate the arithmetic of our maximal quintic $X$.




The case of K3 surfaces shows the existence of quartic surfaces with maximum Picard  number in $\PP^3$.
Explicit models have been derived by Inose in \cite{Inose}. 
In general, surfaces in $\PP^3$ are known to attain the Lefschetz bound only in degree $d\leq 4$ or $d=6$ (see the next section for the latter case). This even holds true if we allow  ADE singularities which is a natural concession since it preserves the deformation type. 

In \cite{Persson}, Persson was able to extend the existence  results for surfaces of maximum Picard  number to certain double covers of rational surfaces. The crucial point about double covers is the following: if the branch curve has at most simple singularities, then the double cover has at most  ADE singularities. Thus one can try to impose enough singularities on the branch curve to obtain a surface with maximum Picard  number as the resolution of the double cover.

Persson mainly considered Horikawa surfaces, i.e.~surfaces attaining Noether's inequality 
\[
K_X^2\geq 2\,p_g(X)-4.
\]
He showed that Horikawa surfaces with maximum Picard  number exist if the congruence condition on the Euler characteristic $\chi\not\equiv 0\mod 6$ is fulfilled.
His approach extends to double covers of $\PP^2$ branched along a curve of arbitrary even degree with at most  simple singularities.

Another construction is due to Bertin and Elencwajg \cite{BE}. For a finite subgroup $G\subset\mbox{Aut}(\PP^1)$, they consider the graphs in $\PP^1\times \PP^1$ of the operation by the group elements. The corresponding conics in $\PP^2$ appear as branch locus of a double cover. This construction gives rise to various projective surfaces of maximum Picard  number.

For elliptic surfaces with section, a uniform picture arises thanks to Shioda's theory of elliptic modular surfaces \cite{Sh-EMS}. 
%
In relation with extremal elliptic surfaces, this approach was generalised by Nori \cite{Nori}.


To our knowledge there is only one other setting where surfaces with maximum Picard  number have turned up so far.
Namely Roulleau studied Fano surfaces parametrising the lines of smooth cubic threefolds.
He derived several instances where the Fano surfaces (which have general type and $h^{2,0}>0$) 
have maximum Picard  number \cite{R1}, \cite{R2}.

It should be pointed out that there are indeed classes of surfaces which do not attain the Lefschetz bound at all. For instance, Livn\'e derived a surface as quotient of the unit ball with $\rho<h^{1,1}$, but without deformations \cite{L}. 

We shall now turn to the quintic surfaces.
The previous record Picard number for quintics with at most ADE singularities was 41 due to Hirzebruch. He considered 5-fold covers of $\PP^1$ branched along five lines. Whenever the intersection points of the lines are distinct, the ten $A_4$ singularities give $\rho\geq 41$ for a minimal desingularisation. 
Actually, Shioda proved in \cite{Sh-PicV} as a consequence of \eqref{eq:Shioda}
that $\rho=41$ for all non-degenerate surfaces in this family.
Thus Theorem \ref{thm} indeed is a genuinely new result.
The next sections elaborate three proofs that $X$ has maximum Picard  number. We shall also investigate the arithmetic of $X$ and determine the zeta function.
In Section \ref{s:rho} we will then consider other Picard numbers of quintic surfaces.

We would like to point out that for numerical quintics (i.e.~smooth minimal surfaces with the same invariants as a smooth quintic in $\PP^3$) Le Barre has constructed an example with maximum Picard number in 1982 \cite{LB}.


\section{Delsarte surfaces}
\label{s:Del}

An irreducible projective surface in $\PP^3$ is called a \emph{Delsarte surface}
if it can be defined by a polynomial which is a sum of four monomials. 
Shioda showed that a Delsarte surface is birational to a Galois quotient $S/G$ of a Fermat surface $S$ by a finite group $G$ \cite{Sh}. He also described an algorithm to find $S$. 
In particular, the transcendental subspace of $H^2(X,\C)$ (the vector subspace generated by transcendental cycles) is identified with the $G$-invariant part of the transcendental subspace of $H^2(S,\C)$.
This enabled Shioda to compute the Picard number $\rho(X)$ in terms of the $G$-action on $S$.

In our case, we can work with the Fermat surface of degree $15$, but we give a general account in terms of the degree $m$:
\[
S_m = \{s^{m}+t^{m}+u^{m}+v^{m}=0\}\subset\PP^3.
\]
The Fermat surface $S_m$ admits coordinate multiplications by $m$-th roots of unity, so projectively $\mu_m^3\subset\mbox{Aut}(S_m)$.
The cohomology of $S_m$ can be decomposed into eigenspaces with character for the induced action of $\mu_m^3$. Here it suffices to consider the following subset of the character group of $\mu_m^3$:
\[
\A_{m}:=
\left\{\alpha=(a_0,a_1,a_2,a_3)\in(\Z/m\,\Z)^4\,|\, a_i\not\equiv 0\pmod{m},\,\sum_{i=0}^3 a_i\equiv 0\pmod{m}\,\right\}.
\]
For $\alpha\in\A_{m}$, let $V(\alpha)$ denote the corresponding eigenspace with character. Here we let $g=(\zeta_1, \zeta_2, \zeta_3)\in\mu_m^3$ operate on $S_m$ as
\begin{eqnarray}
\label{eq:op}
 [s,t,u,v] \mapsto [s,\zeta_1\,t, \zeta_2\,u,\zeta_3\,v].
\end{eqnarray}
Then the subspace $V(\alpha)\subset H^2(S_m)$ is determined by the condition
\[
g^*|_{V(\alpha)} = \alpha(g) = \zeta_1^{a_1}\,\zeta_2^{a_2}\,\zeta_3^{a_3} \;\;\; \forall\; g=(\zeta_1, \zeta_2, \zeta_3)\in\mu_m^3.
\]

By results of Katz and Ogus, each $V(\alpha)$ is one-dimensional 
(this holds true for Fermat varieties of arbitrary dimension).
One has
\begin{eqnarray}
\label{eq:esp}
H^2(S_m) = V_0 \oplus \bigoplus_{\alpha\in\A_{m}} V(\alpha)
\end{eqnarray}
where $V_0$ corresponds to the trivial character and is spanned by the hyperplane section.

We briefly explain how to decide whether $V(\alpha)$ is algebraic or transcendental.
Consider the subspace $H^{2,0}(S_m)\subset H^2(S_m,\C)$.
In the affine chart $s=1$, $H^{2,0}(S_m)$ is generated by the following $2$-forms:
\[
\omega(b_1,b_2,b_3)=t^{b_1-1}u^{b_2-1}v^{b_3-1}\frac{du\wedge dv}{t^{m-1}},\;\; b_i\geq 1,\;  b_1+b_2+b_3\leq m-1.
\]
An automorphism $g=(\zeta_1, \zeta_2, \zeta_3)\in\mu_m^3$ acts on these $2$-forms by
\[
g^*\omega(b_1,b_2,b_3) = \zeta_1^{b_1}\,\zeta_2^{b_2}\,\zeta_3^{b_3} \omega(b_1,b_2,b_3).
\]
Let $b_0=m-(b_1+b_2+b_3)$ and $\alpha=(b_0,b_1,b_2,b_3)$ 
where we abuse notation by not distinguishing between the integers $b_i$ with $0<b_i<m$ and their equivalence classes in $\Z/m\Z$.
The eigenspace decomposition \eqref{eq:esp} implies that $V(\alpha)=\C\omega(b_1,b_2,b_3)$.
It follows that
\[
H^{2,0}(S_m)\oplus H^{0,2}(S_m) = \bigoplus_{\alpha\in\TT_m} V(\alpha) \subset H^2(S_m,\C)
\]
where
\[
\TT_m = \left\{\alpha=(b_0,b_1,b_2,b_3)\in\A_m; \; 0<b_i<m,\; \sum_{i=0}^3 b_i = m \text{  or } 3m\right\}
\]
The eigenspace decomposition \eqref{eq:esp} is defined over $\Q(\zeta_m)$.
Here the Galois group $\mathfrak{G}=\mbox{Gal}(\Q(\zeta_m)/\Q)\cong(\Z/m\,\Z)^*$ 
operates on $\A_{m}$ coordinatewise by multiplication. 
The space of transcendental cycles $T(S_m)\subset H^2(S,\Q)$ 
is the smallest $\Q$-vector subspace $V$
such that $V\otimes \C$ contains $H^{2,0}(S_m)\oplus H^{0,2}(S_m)$.
Thus we find
\[
T(S_m) = \bigoplus_{\alpha\in\mathfrak{G}\TT_{m}} V(\alpha).
\]

\begin{Example}[Fermat Quintic]
\label{Ex:Fermat}
{\rm
A classical example is the Fermat quintic $S_5$.
One easily finds that $\mathfrak{G}\TT_5$ consists of four $(\Z/5\Z)^*$ orbits corresponding to the element $(1,1,1,2)\in\TT_5$ and the coordinate permutations. Hence $\dim(T(S_5))=16$ and $\rho(S_5)=37$. Since $h^{1,1}(S_5)=45$ as in the introduction, $S_5$ does not have maximum Picard  number.}
\end{Example}

One can easily show that in higher degree $m>5$, the Fermat surface $S_m$ has maximum Picard  number if and only if $m=6$. In fact, the $(\Z/m\,\Z)^*$-orbit of $(1,1,1,m-3)\in\TT_m$ contains a character with eigenspace of Hodge weight $(1,1)$ if and only if $\phi(m)>2$. By definition, this eigenspace is non-algebraic for $m>3$. 
Alternatively, one can compare the asymptotic growth of $\rho(S_m)$ as $3m^2$ (cf.~\cite{Sh-Pic}) against $h^{1,1}$ which is asymptotic to $2m^3/3$. The exceptional property of the Fermat sextic was noticed by Beauville.

By definition, a Delsarte surface is covered by a suitable Fermat surface. Shioda gave an algorithm to find the Fermat degree $m$ and the dominant rational map $\varphi$ \cite{Sh}. 
In case of the quintic  $X$ from Theorem \ref{thm}, one finds $m=15$ and 
\begin{eqnarray*}
\varphi:\;\;\;\;\;\;  S_{15}\;\;\;\;\;\; & \dasharrow & \;\;\;\;\;\;\; X\\
{[s,t,u,v]} & \mapsto & [t\,u^3\,v^7, s\,t^3\,u^7, v\,s^3\,t^7, u\,v^3\,s^7].
\end{eqnarray*}

The Delsarte surface $X$ is birational to the quotient $S_m/G$ where $G$ is the covering group corresponding to $\varphi$, i.e.~$G=\{g\in\mu_m^3; \varphi=\varphi\circ g\}$. Since the Lefschetz number
\[
\lambda(X) = b_2(X) - \rho(X)
\]
is a birational invariant, we can compute it (and thus $\rho(X)$) through the quotient $S_m/G$. Let $\TT_{m}^G$ consist of all those $\alpha\in\TT_{m}$ such that all elements in $G$ act as identity on
$V(\alpha)$. This is computed as follows: Write $G\ni g=(\zeta_1, \zeta_2, \zeta_3)$, operating on $S_m$ as in \eqref{eq:op}.
Let $\alpha=(a_0,a_1,a_2,a_3)\in\mathfrak{A}_m$. Then $V(\alpha)$ is $G$-invariant if and only if 
\[
\prod_{i=1}^3 \zeta_i^{a_i} =1 \;\;\;\;\;\forall\, g=(\zeta_1, \zeta_2, \zeta_3)\in G.
\]

For the Lefschetz number, we obtain
\[
\lambda(X)=\lambda(S_m/G) = \#\mathfrak{G}\TT_{m}^G.
\]
In our case, one easily finds that $\mathfrak{G}\TT_{15}^G$ is the $(\Z/15\,\Z)^*$ orbit of a single element, say $\alpha=(1,2,4,8)$. Hence $\lambda(X)=8$ and $\rho(X)=45$ as claimed in Thm.~\ref{thm}. \qed

\section{Generators of the N\'eron-Severi group}

In this section, we work out an explicit $\Q$-basis of the N\'eron-Severi group of our quintic surface $X$.
This gives an alternative proof of Theorem \ref{thm}
and  enables us to compute the zeta function of $X$ in the next section.

We first have to consider the resolution of singularities on $Y$. It is easily checked that the only singularities occur at $[0,0,0,1]$ and permutations, and that they have type $A_9$. Hence we already have $\rho(X)\geq 37$. 

We consider three further groups of rational curves on $X$:

\begin{enumerate}
 \item 
The strict transforms of the six lines in $\PP^3$ passing through any two nodes of $Y$:
\[
\ell_{xy}=\{x=y=0\}\subset\PP^3,\;\;\; \ell_{xz}=\hdots.
\]

\item
The five lines
\[
 \ell_\alpha = \{x=\alpha\, z, y=\alpha^7\, w\}\subset X,\;\;\;\; \alpha^5=-1.
\]

\item
The images of the non-contracted lines on $S_{15}$
\begin{eqnarray*}
 C_\varrho=\{[\varrho^i\,\mu^3,-\lambda\,\mu^2, \varrho^i\,\lambda^3,-\mu\,\lambda^2]; \; [\lambda,\mu]\in\PP^1\}, &\\
D_\varrho=\{[-\lambda\,\mu^2, \varrho^i\,\lambda^3,-\mu\,\lambda^2,\varrho^i\,\mu^3]; \; [\lambda,\mu]\in\PP^1\}, &\;\;\varrho^3=1.
\end{eqnarray*}
\end{enumerate}

The intersection behaviour with the exceptional locus is sketched in the following figure for the node $[0,0,0,1]$. 
Here we number the components of the exceptional divisor from $1$ to $9$ while $D_\varrho$ stands for all three rational curves with $\varrho^3=1$.

$$
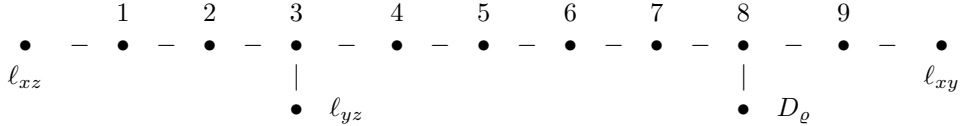
\begin{array}{ccccccccccccccccccccc}
&& 1 && 2 && 3 && 4 && 5 && 6 && 7 && 8 && 9 &&\\
\bullet & - & \bullet & - & \bullet & - & \bullet & - & \bullet & - & \bullet & - &
\bullet & - & \bullet & - & \bullet & - & \bullet & - & \bullet\\
\ell_{xz} &&&&&& | &&&&&&&&&& | &&&& \ell_{xy}\\
&&&&&& \bullet & \ell_{yz} &&&&&&&&& \bullet & D_\varrho & &&\\
\end{array}
$$

The verification is straight forward by computing the resolution of the $A_9$ singularity.
The intersection behaviour at the other nodes is obtained by cyclic permutation of coordinates
\[
[x,y,z,w] \mapsto [w,x,y,z].
\]
All other non-zero intersection numbers are given as follows:
\[
 C_\varrho.D_{\varrho^2}=5, \;\;\; C_\varrho.\ell_\alpha=D_\varrho.\ell_\alpha=1,\;\; C_\varrho.\ell_{xz}=C_\varrho.\ell_{yw}= D_\varrho.\ell_{xz}=D_\varrho.\ell_{yw}=1.
\]

Finally for the self-intersection numbers, we let $H$ denote the hyperplane section. Then $\ell_*.H=1, C_\varrho.H=D_\varrho.H=3$. Hence the adjunction formula with $K_X=H$ gives
\[
\ell_*^2=-3,\;\; C_\varrho^2=D_\varrho^2=-5.
\]

We will now exhibit a $\Q$-basis of $\NS(X)$. Consider the following 45 rational curves on $X$:
\[
\mathcal B = \{4 \times A_9,\, \ell_{xy},\, \ell_{yz},\, \ell_{xz},\, C_\varrho\, (\varrho\neq 1),\, \ell_\alpha\, (\alpha\neq -1)\}.
\]
Their  intersection matrix has determinant $202500=2^2\, 3^4\, 5^4$. Since $\rho(X)\leq 45$ by Lefschetz' bound (\ref{eq:Lef}), we deduce $\rho(X)= 45$. The above curves give a $\Q$-basis of $\NS(X)$, i.e.~they generate $\NS(X)$ up to finite index. \qed

\begin{Remark}
{\rm
A joint paper with Shioda and van Luijk introduced a supersingular reduction technique to prove that $\NS(S_m)$ is integrally generated by lines for all $m\leq 100$ that are relatively prime to $6$ \cite{SS}. The same method is applicable here for $X$. One could try to work with the supersingular reduction at $p=29$.}
\end{Remark}

\section{Zeta function}
\label{s:zeta}

We are now in the position to determine the zeta function of $X$. We will deal with the algebraic part $\NS(X)$ and the transcendental part $T(X)$ separately.

For the algebraic part, we consider $\NS(X)$ as a subspace of $H^2(X)$ in some \'etale cohomology. Hence the eigenvalues of Frobenius are $p$ times roots of unity.
Note that the rational basis $\mathcal B$ is Galois invariant. Hence the contribution of $\NS(X)$ to the zeta function is as follows:

\begin{Lemma}\label{Lem:zeta-NS}
Let $K$ resp.~$L$ denote the third resp.~fifth cyclotomic field over $\Q$. Then
\[
L(\NS(X), s) = \zeta_{\Q}(s-1)^{39} \,\zeta_K(s-1)\,\zeta_L(s-1).
\]
\end{Lemma}

For the transcendental part, Weil translated the motivic decomposition of $H^2(S_m)$ into Jacobi sums \cite{Weil}. We follow his description of the local Euler factors for a suitable prime power $q=p^r$ such that
\[
q\equiv 1 \mod m.
\]
On the field $\F_q$ of $q$ elements, we fix a character 
\[
\chi: \F_q^* \to \C^*
\]
of order exactly $m$. For any $\alpha\in\mathfrak{A}_m$, we then define the Jacobi sum
\begin{eqnarray}\label{s:Jacobi}
j(\alpha) = \sum_\text{\small $\begin{matrix} v_1, v_2, v_3\in\F_q^*\\ v_1+v_2+v_3=-1\end{matrix}$} \chi(v_1)^{a_1} \chi(v_2)^{a_2}\chi(v_3)^{a_3}.
\end{eqnarray}

\begin{Theorem}[Weil]
\label{Thm:Weil}
In the above notation, consider the Fermat surface $S_m$ over $\F_q$ with Frobenius morphism Frob$_q$. Then Frob$_q^*$ has the following characteristic polynomial on $H^2(S_m)$:
\[
P(T) = (T-q)  \prod_{\alpha\in\mathfrak{A}_m} (T-j(\alpha)).
\]
\end{Theorem}

We will now use Theorem \ref{Thm:Weil} to determine the local Euler factors of the transcendental subspace $T(X)$. We are concerned with the covering Fermat surface $S_{15}$. By section \ref{s:Del}, $T(X)$ is identified with a single $(\Z/15\Z)^*$-orbit
\[
 T(X) = \bigoplus_{\alpha\in\mathfrak{G}\TT_{15}^G} V(\alpha) = \bigoplus_{k\in(\Z/15\Z)^*} V(k\cdot(1,2,4,8)).
\]

Since the dominant rational map $S_m\to X$ is defined over $\Q$, we obtain

\begin{Lemma}\label{Lem:zeta-T}
Let $q\equiv 1\mod 15$. Then the local Euler factor of $T(X)$ at $q$ is
\[
 L_q(T(X), s) =  \prod_{\alpha\in\mathfrak{G}\TT_{15}^G} (1-j(\alpha)\, q^{-s}).
\]
\end{Lemma}

Together, Lemma \ref{Lem:zeta-NS} and \ref{Lem:zeta-T} determine the zeta function of $X$:

\begin{Proposition}
Let $L(T(X), s)$ denote the $L$-series of $T(X)$ as given by the local Euler factors in Lemma \ref{Lem:zeta-T}. Then 
\[
 \zeta(X, s) = \zeta_\Q(s)\,\zeta_{\Q}(s-1)^{39} \,\zeta_K(s-1)\,\zeta_L(s-1)\,L(T(X),s)\,\zeta_\Q(s-2).
\]
\end{Proposition}

\section{Automorphisms}

The third proof of Theorem \ref{thm} could be considered most ad hoc,
as it requires the least information about the surface $X$.
The basic idea is to combine the existence of an automorphism of order $15$ on $X$ (which comes of course from the covering Fermat surface $S_{15}$) with just a little knowledge about $\NS(X)$. Here the operation of the automorphism on the holomorphic 2-forms on $X$ will enable us to see $\rho(X)=45$ easily.

The quintic surface $X$ admits an automorphism $g$ of order $15$. Let $\zeta$ denote a primitive $15$th root of unity. Then $g$ can be given by
\[
g(x,y,z,w) = [\zeta\,x, \zeta^3\,y, \zeta^7\,z, w]
\]
We determine the operation of $g$ on $H^{2,0}(X)$. We express a basis of $H^{2,0}(X)$ in the affine chart $w=1$ in terms of
\[
\omega = \dfrac{dy\wedge dz}{\partial_x F} = \dfrac{dy\wedge dz}{y\,z^3+y^3+3\,z\,x^2}.
\]
By Griffiths' residue theorem, a basis of $H^{2,0}(X)$ and the operation of $g^*$ is as follows:
$$
\begin{array}{c||c|c|c|c}
\hline
\text{basis} & \omega & x\,\omega & y\,\omega & z\,\omega\\
\hline
g^* & \zeta & \zeta^2 & \zeta^4 & \zeta^8\\
\hline
\end{array}
$$
For our purposes, it is crucial that these eigenvalues amount for exactly half of all complex embeddings $\Q(\zeta)\hookrightarrow{\C}$.
Since there are no conjugate duplicates involved, the eigenvalues in fact form a CM-type of $\Q(\zeta)$. 
It follows that $g^*$ endows $T(X)$ with the structure of a $\Q[\zeta]$-vector space.
In particular 
\begin{eqnarray}\label{eq:phi}
 8=\phi(15)\mid\mbox{dim}(T(X)).
\end{eqnarray}
Here the four $A_9$ singularities on $Y$ give $\rho(X)\geq 37$, so $T(X)$ has dimension $8$ or $16$. 
In fact, taking the strict transforms of any two distinct lines through two nodes of $Y$, we see $\rho(X)\geq 38$ and $\dim(T(X))\leq 15$. By (\ref{eq:phi}), this implies dim$(T(X))=8$ and thus $\rho(X)=45$. This completes the third proof of Theorem \ref{thm}. \qed

\begin{Remark}
The ideas from this section 
can be employed to search for surfaces in $\PP^3$ with maximum Picard number 
in a systematic manner. 
However, for degree $d>4$,
we did not find any surfaces with only ADE-singularities other than $X$ up to isomorphism.
\end{Remark}

\section{Smaller Picard numbers}
\label{s:rho}

We will now consider quintic surfaces with smaller Picard numbers.
Some examples were given by Shioda in \cite{Sh-PicV}.
Note that all those Picard numbers are congruent to $1$ modulo $4$.
Here we shall exhibit quintic surfaces with several further Picard numbers.

We employ a systematic approach through Delsarte surfaces.
Namely we isolate all quintic Delsarte surfaces with only ADE-singularities.
Then we compute their Picard numbers using the technique from Section \ref{s:Del}.
Notably we will also find odd Picard numbers congruent to $3$ modulo $4$ (as indicated in Theorem \ref{thm2}).

To exclude the Delsarte surfaces with singularities worse than  rational double points
we proceed as follows.
We have already pointed out that a smooth quintic $X$ or the minimal desingularisation of a quintic with only  rational double points has $h^{2,0}(X)=4$.
If there are worse singularities, 
then this necessarily causes $h^{2,0}$ to drop.
We exclude those quintic Delsate surfaces by considering the $G$-invariant eigenspaces $V(\alpha)$ on the covering Fermat surface $S_m$. 
As explained in Section \ref{s:Del},
the Hodge type of the eigenspace $V(\alpha)$ is determined by the reduced representative $\alpha=(b_0,\hdots,b_3)$ with $0<b_i<m$ in terms of 
\[
|\alpha| = (b_0+\hdots+b_3)/m-1.
\]
Namely $V(\alpha)$ has Hodge type $(2-|\alpha|, |\alpha|)$.
For a quintic Delsarte surface, we thus find the invariant eigenspaces $V(\alpha)$ of Hodge type $(2,0)$
on the covering Fermat surface,
and we can check whether there are exactly four of them.

The next table collects all Picard numbers that arise from quintic Delsarte surfaces with  rational double points.
For each, we give a defining polynomial for a  quintic surface with this Picard number.
In the known cases, the last column refers to \cite{Sh-PicV}, although in two cases ($\rho=17,41$) we decided to include explicit new examples as opposed to the generic examples in \cite{Sh-PicV}.
In the new cases, the last column of the table specifies the ADE-types of the singularities.


\begin{table}[ht!]
$$
\small{
\begin{array}{ccc}
\hline
\text{Picard number} & \text{polynomial} & \text{comment}\\
\hline
\rho=1 & 
xy^4+yz^4+zx^4+w^5  & \cite[\text{Thm.}~4.1]{Sh-PicV}\\
\rho=5 & 
x^5+xy^4+yz^4+w^5 & \cite{Sh-PicV}\\
\rho=13 & x^5+y^5+xzw^3+wz^4 & A_4\\
\rho=17 & wx^4+wy^4+yz^4+zw^4 & 4A_3\\
\rho=19 & ywx^3+xy^4+yz^4+zw^4 & A_{16}\\
\rho=21 & 
xy^4+yz^4+zw^4+wx^4 & \cite{Sh-PicV}\\
\rho=23 &
ywx^3+y^5+wz^4+zw^4 & A_{19}\\
\rho=25 & 
x^5+xy^4+z^5+w^5 & \cite{Sh-PicV}\\
\rho=27 & yzx^3+wy^4+z^5+w^5 & A_4\\
\rho=29 & 
x^5+xy^4+z^5+zw^4 & \cite{Sh-PicV}\\
\rho=31 & zw^4+yz^4+xzy^3+ywx^3 & A_{12}+A_{16}\\
\rho=33 & ywx^3+zwy^3+yz^4+w^5 & A_{11}+A_{19}\\
\rho=35 & ywx^3+wy^4+wz^4+zw^4 & 4A_2+A_{16}\\
\rho=37 & x^5+y^5+z^5+w^5 & \text{Ex.}~\ref{Ex:Fermat}\\
\rho=39 & yzx^3+wy^4+wz^4+w^5 & 4A_2+A_4\\
\rho=41 & 
xy^4+xz^4+zx^4+zw^4 & 8A_3\\
\rho=43 & 
zw^4+wz^4+wzy^3+yx^4 & 7A_3\\
\rho=45 & yzw^3+xyz^3+wxy^3+zwx^3 & \text{Thm.}~\ref{thm}\\
\hline
\end{array}}
$$
\caption{Quintic surfaces and their Picard numbers (after desingularisation)}
\end{table}

One can check that a quintic Delsarte surface with Picard number $\rho=45$
is unique up to trivial coordinate change,
provided its singularities are only  rational double points.
Such a  uniqueness result does not hold for quintic Delsarte surfaces with smaller Picard number 
(see e.g.~\cite{Sh-PicV}).

There are five small odd Picard numbers missing in the table
(as specified in Theorem \ref{thm2}) 
as well as all even Picard numbers.
To overcome this lack of explicit examples, we have recently started a project with R.~van Luijk where we aim at engineering quintic surfaces with prescribed Picard number explicitly.


\subsection*{Acknowledgements}
The problem of finding complex surfaces with maximum Picard  number was brought to our attention by T.~Shioda.
We are grateful to him and to B.~van Geemen, X.~Roulleau and R.~van Luijk for many helpful comments and discussions.
Our thanks go to the referee for many suggestions that helped improve the paper.
When this paper was started, the author held a position at the Mathematics Department of the University of Copenhagen.
 

\vspace{0.8cm}

Matthias Sch\"utt\\
Institut f\"ur Algebraische Geometrie\\
Leibniz Universit\"at Hannover\\
Welfengarten 1\\
30167 Hannover\\
Germany\\
schuett@math.uni-hannover.de

\end{document}